\documentclass{conm-p-l}

\newtheorem{theorem}{Theorem}[section]

\newtheorem{prop}[theorem]{Proposition}

\theoremstyle{definition}
\newtheorem{definition}[theorem]{Definition}
\newtheorem{example}[theorem]{Example}

\theoremstyle{remark}

\numberwithin{equation}{section}

\usepackage{amsmath}
\usepackage{amssymb}

\renewcommand{\bar}{\overline}

\newcommand{\CC}{\mathbb{C}}

\newcommand{\FF}{\mathbb{F}}

\newcommand{\QQ}{\mathbb{Q}}
\newcommand{\RR}{\mathbb{R}}

\newcommand{\TT}{\mathbb{T}}

\newcommand{\ZZ}{\mathbb{Z}}

\newcommand{\Qp}{\QQ_p}
\newcommand{\Zp}{\ZZ_p}

\newcommand{\Fp}{\FF_p}

\newcommand{\ints}{{\mathcal O}}

\newcommand{\calC}{{\mathcal C}}
\newcommand{\calD}{{\mathcal D}}
\newcommand{\calE}{{\mathcal E}}

\DeclareMathOperator{\Tr}{Tr}
\DeclareMathOperator{\Supp}{Supp}

\begin{document}


\title{Examples of Wavelets for Local Fields}
\author{Robert L. Benedetto}
\address{Department of Mathematics and Computer Science \\
	Amherst College \\
	Amherst, MA 01002 \\
	USA}
\email{rlb@cs.amherst.edu}
\urladdr{http://www.cs.amherst.edu/\textasciitilde rlb}

\date{June 18, 2003}



\begin{abstract}
It is well known that the Haar and Shannon wavelets in $L^2(\RR)$
are at opposite extremes, in the sense that the Haar wavelet is
localized in time but not in frequency, whereas the Shannon wavelet
is localized in freqency but not in time.  We present a
rich setting where the Haar and Shannon wavelets
coincide and are localized both in time and in frequency.

More generally, if $\RR$ is replaced by a group $G$ with
certain properties, J.~Benedetto and the author have
proposed a theory of wavelets on $G$, including the construction
of wavelet sets \cite{BB}.  Examples of such groups $G$
include the $p$-adic rational
group $G=\Qp$, which is simply the completion of $\QQ$ with respect
to a certain natural metric topology, and the Cantor dyadic
group $\FF_2((t))$ of formal Laurent series with coefficients
0 or 1.

In this expository paper, we consider
some specific examples of the wavelet theory on such groups $G$.
In particular,
we show that Shannon wavelets on $G$ are the same as Haar wavelets
on $G$.  We also give several examples of specific groups (such
as $\Qp$ and $\FF_p((t))$, for any prime number $p$)
and of various wavelets on those groups.
All of our wavelets are localized in frequency; the Haar/Shannon
wavelets are localized both in time and in frequency.
\end{abstract}

\maketitle

One of the principal goals of wavelet theory has been the
construction of useful orthonormal bases for $L^2(\RR^d)$.  The group
$\RR^d$ has been an appropriate setting both because of its use in
applications and because of its special property of containing
lattices such as $\ZZ^d$ which induce discrete groups of translation
operators on $L^2(\RR^d)$.
Of such wavelets, the easiest to describe are the Haar wavelet
and the Shannon wavelet in $L^2(\RR)$.  The Haar is a compactly
supported step function in time but has noncompact support and
slow decay in frequency.  On the other hand, the Shannon is
considered to be at the opposite extreme, being a compactly
supported step function in frequency but with noncompact support
and slow decay in time.

This distinction between Haar and Shannon wavelets appears to
caused by certain properties of the topological group $\RR^d$.
J.~Benedetto and the author \cite{BB} have presented a theory
of wavelets on $L^2(G)$, for groups $G$ with different properties
to be described in Section~\ref{sect:group} below.
In this paper, we shall present some of
that theory and some examples.  Our main result,
Theorem~\ref{thm:haarshan}, is the surprising
fact that Haar and Shannon wavelets not only
exist in $L^2(G)$, but the two are the same.  Moreover, these
Haar/Shannon wavelets are localized both in time and in frequency;
in fact, both the wavelets and their transforms are compactly
supported step functions.

To reach that goal,
we must first generalize wavelet theory for $\RR^d$
to a more abstract setting.
A number of authors have
extended wavelet theories to locally compact groups $G$
other than $\RR^d$;
see, for example, \cite{Dah,Far,HLPS,Hol,Luk,PS}.
However, those generalizations required the group $G$
to contain a discrete subgroup, just as $\RR^d$ contains
the discrete lattice $\ZZ^d$.
The more abstract generalized multiresolution analyses in
\cite{BMM,Pap} are even broader, though the discrete group of
translation operators is, naturally, a requirement and indeed
one of the defining characteristics of a wavelet theory.
Indeed, the standard constructions of multiresolution
analysis and of wavelet sets in $L^2(\RR^d)$
rely crucially on a discrete group of translation
operators (as can be seen in expositions and investigations
such as \cite{BL1,BL2,DLS,Dau,Mey,SW}).

Beyond the identification of Haar and Shannon wavelets, we
are also motivated by the study of analysis on local fields,
which are of great importance in number theory; see
\cite{Gou,Kob,Rob,Ser} for background.  The study
of Fourier analysis on local fields and on ad\`{e}le groups
arises in class field theory and in the study of certain
zeta functions; see \cite{Con1,Con2,RV,Tate}, for example.
However, no prior knowledge of such topics are needed to read
this paper.

One simple example of a local field we wish to consider
is the field $\FF_p((t))$
of formal Laurent series over the field $\FF_p$ of
$p$ elements, where $p$ is a prime;
see Examples~\ref{ex:FpT},~\ref{ex:FpTD},
and~\ref{ex:FpTwave}.  As an additive group, $\FF_p((t))$
contains a compact open subgroup $\FF_p[[t]]$ (consisting
of formal power series) as well as a discrete subgroup
$\Gamma$, with $\FF_p((t)) = \FF_p[[t]] + \Gamma$.
Using the elements of $\Gamma$ as translations, one
can study wavelets in $L^2\left(\FF_p((t))\right)$.  For example,
Lang \cite{Lan} examined the MRA theory in the case
$p=2$; we shall consider this case briefly in
Section~\ref{sect:motiv} to motivate the generalizations
to follow.

However, while some of the groups $G$ we wish to consider
contain a discrete lattice, others do not.  For example,
the $p$-adic field $\Qp$, which we shall study in
Examples~\ref{ex:Qp},~\ref{ex:QpD},~\ref{ex:Qpwave},
and~\ref{ex:Qpwave3}, is a completion of $\QQ$ with respect
to a certain metric and is a natural analogue of the
real line $\RR$.  Unfortunately, as an additive
group, $\Qp$ contains no nontrivial discrete subgroup, thus
presenting a serious obstacle to the development of a theory
of wavelets in $L^2(\Qp)$.  Kozyrev \cite{Koz} produced one
specific set of wavelet generators, analogous to Haar wavelets
on $\RR^d$, using a discrete set of translation operators
which do not form a group.  However, as we shall see, those
operators allow the Haar wavelets, but they
preclude the possibility of a general theory
of wavelet sets or multiresolution analysis in $L^2(\Qp)$.

Fortunately, although it lacks a lattice, $\Qp$ does contain
the compact open subgroup $\Zp$.  Whereas $\RR^d$ contains
a discrete lattice $\ZZ^d$ with compact quotient $\RR^d/\ZZ^d$,
then, the group $\Qp$ contains a {\em compact} subgroup $\Zp$
with {\em discrete} quotient $\Qp/\Zp$.  Thus, the situation
for $\Qp$ is the reverse of that for $\RR^d$, but it still
permits the construction of a discrete {\em group} of operators
on $L^2(\Qp)$.  The resulting operators $\tau_{[s]}$
are not simple translation operators, but one may compute
with them using the results of Section~\ref{sect:oper} below.

Of course, one also needs a dilation, which we
shall call an expansive automorphism, to be defined
in Section~\ref{sect:oper}.
J.~Benedetto and the author \cite{BB} therefore
considered wavelets on
any locally compact abelian group $G$ with
a compact open subgroup $H$ and an expansive automorphism $A$.
The resulting wavelets, developed independently
of Kozyrev's work, form a large class of orthonormal
wavelet bases over $\Qp$, $\FF_p((t))$, and any other
such group $G$.  The theory also unifies both Lang's
and Kozyrev's wavelets within a much broader setting.

Besides identifying Haar and Shannon wavelets in
this context,
the main purpose of this
paper is to give examples of some of
the groups and constructions considered in \cite{BB} in
an expository setting.
In Section~\ref{sect:motiv}, we present a concrete
example of a Haar/Shannon wavelet, in the Cantor dyadic
group $\FF_2((t))$.
In Section~\ref{sect:group}, we list some basic properties
of a group $G$ with compact open subgroup $H$ and
expansive automorphism $A$; we also present several
examples of such groups.
In Section~\ref{sect:oper} we describe the
discrete group of translation-like operators on $L^2(G)$
defined in \cite{BB},
and we consider specific examples by revisiting several of
the groups from Section~\ref{sect:group}.
In Section~\ref{sect:haar} we define wavelets using
those operators as in \cite{BB}.
We also present an explicit analogue of Shannon
wavelets on $G$ in Theorem~\ref{thm:haardef}.  Then,
in Theorem~\ref{thm:haarshan},
we prove the promised result that our
Shannon wavelets {\em are} Haar wavelets.
Finally, in Section~\ref{sect:exs}, we present several
examples of wavelet sets produced by the results of \cite{BB}
over certain groups $G$.

\section{A motivating example}
\label{sect:motiv}
We begin by observing the Haar/Shannon phenomenon in
a specific and relatively simple context.
Let $\FF_2=\{0,1\}$ be the field of two elements, with
operations $0+0=1+1=0$, $1+0=0+1=1$,
$0\cdot 0 = 0\cdot 1 = 1\cdot 0 = 0$, and
$1\cdot 1 = 1$.  The field $G=\FF_2((t))$ of formal Laurent
series consists of all infinite formal sums
$$c_{n_0} t^{n_0} + c_{n_0 + 1} t^{n_0 + 1} 
+ c_{n_0 + 2} t^{n_0 + 2} + \cdots,$$
where $n_0\in\ZZ$, each $c_n\in\FF_2$, and no consideration
is given to convergence or divergence.  One can add or
multiply elements of $G$ as formal power series: addition
follows the rule $b_n t^n + c_n t^n = (b_n + c_n)t^n$,
and multiplication follows the rule
$(c_n t^n) \cdot (c_m t^m) = (c_n \cdot c_m) t^{n+m}$.
One can also put a metric topology on $G$ by declaring
the absolute value of $g\in G$ to be $|g|=2^{-n_0}$,
where $n_0$ is the smallest integer such that $c_{n_0}\neq 0$.
(We will see the more general fields $\FF_p((t))$ later
in Examples~\ref{ex:FpT},~\ref{ex:FpTD}, and~\ref{ex:FpTwave}.)

Lang \cite{Lan} considers the ``Cantor dyadic group'',
which is the group $G=\FF_2((t))$ under addition.  He
views $G$ as the subgroup of $\prod_{n\in\ZZ} \FF_2$
consisting of elements $(\cdots,c_{-2},c_{-1},c_0,c_1,c_2,\ldots)$
for which there exists an $n_0\in\ZZ$ such that
$c_n=0$ for all $n<n_0$.  The metric topology we have described
on $G$ coincides with the product topology $G$ inherits
from $\prod_{n\in\ZZ} \FF_2$.

One can show that $G$ is self-dual.  Whereas the duality
pairing on $\RR$ takes $x,y\in \RR$ to $e^{2\pi i xy}$,
the pairing on $G$ takes $g_1,g_2\in G$ to $e^{\pi i c}$,
where $c=c_{-1}$ is the $t^{-1}$-coefficient of $g_1\cdot g_2$.
Thus, the pairing of $g_1$ and $g_2$ is $1$ if $c=0$,
and $-1$ if $c=1$.

$G$ contains a discrete subgroup $\Gamma$ consisting of elements
for which $c_0=c_1=c_2=\cdots = 0$.  (That is, $\Gamma$ consists
of all the formal sums which have only finitely many nonzero
coefficients, and all nonnegative power terms $c_n t^n$ are zero.)
The compact subgroup $H=\FF_2[[t]]$ of $G$, consisting
of all power series $c_0 + c_1 t + c_2 t^2 + \cdots$, is
a fundamental domain for $\Gamma$, in the sense that
$H\cap\Gamma=\{0\}$ and $H+\Gamma = G$.  It may be helpful
to think of $\Gamma$ as analogous to the discrete subgroup
$\ZZ$ of $\RR$, while $H$ is analogous to the fundamental
domain $[-1/2,1/2)$; however, in our case, $H$ is itself
a subgroup of $G$.

We may normalize a Haar measure $\mu$ on $G$ with $\mu(H)=1$.
Note that $tH=\{c_1 t + c_2 t^2 + \cdots\}$ is a subgroup of $H$
of measure $1/2$, because $H$ is the disjoint union
$(1+tH) \cup tH$, and $\mu(tH)=\mu(1+tH)$ by translation invariance.
Note also that $H$ is homeomorphic to the Cantor set; it consists
of two disjoint pieces, $tH$ and $1+tH$, each of which is homeomorphic
(via multiplication-by-$t^{-1}$) to $H$.

Similarly, $t^{-1}H=\{c_{-1} t^{-1} + c_0 + c_1 t + \cdots\}$
is a subgroup of $G$ which contains $H$ and which has measure $2$.
It is easy to see that
$$\cdots \supset t^{-2}H \supset t^{-1} H \supset H
\supset tH \supset t^2 H \supset \cdots,$$
and
$$\bigcap_{n\geq 0} t^n H = \emptyset, \quad\text{and} \quad
\bigcup_{n\geq 0} t^{-n} H = G.$$
Thus, multiplication-by-$t^{-1}$ is a dilation which expands areas
by a factor of $2$, just as multiplication-by-$2$ is a dilation
of $\RR$ which expands areas by a factor of $2$.

Lang \cite[Example 5.2]{Lan} observed that the function
$$f(g)=
\begin{cases}
1, & \text{if } g\in tH ; \\
-1, & \text{if } g\in 1+tH ; \\
0, & \text{otherwise}
\end{cases}
$$
is a wavelet, in the sense that the set of dilated translates
$\{f(t^n g - s) : n\in\ZZ, s\in\Gamma\}$ is an orthonormal
basis for $L^2(G)$.  Indeed, this function $f$ is an obvious
analogue of the usual Haar wavelet on $\RR$.  The surprising
fact that motivates much of the rest of this paper is
that using the duality pairing on $G$, one can easily
compute that the Fourier transform of $f$ is
$$\widehat{f}(\omega)=
\begin{cases}
1, & \text{if } \omega \in t^{-1} + H ; \\
0, & \text{otherwise,}
\end{cases}
$$
which shows that $\widehat{f}$ is a step function with compact support
in the frequency domain, and that it is therefore
natural to call $f$ a Shannon wavelet.
(The skeptical reader may observe that $t^{-1}+H$ may be
written as the disjoint union of $t^{-1} + (tH)$ and
$-t^{-1} + (1+tH)$, just as the Shannon wavelet in $L^2(\RR)$
has Fourier transform equal to the characteristic function
of $[-1,-1/2)\cup[1/2,1)$.)

It turns out that the key property of $G$ that allows the
Haar/Shannon identification is the fact that the fundamental
domain $H$ is itself a subgroup of $G$.  Moreover, the
lattice $\Gamma$ turns out to be unnecessary, provided we
define operators to act as translations in some appropriate
manner.  Thus, we will be able to work with groups such
as $\Qp$, which have no such lattice.
With these ideas in mind, we are now prepared
to generalize.

\section{The group, subgroup, and automorphism}
\label{sect:group}
We fix the following notation.
\begin{tabbing}
$G$ \= \hspace{1.0in} \= a locally compact abelian group \\
$\widehat{G}$ \> \> the dual group of $G$ \\
$H$ \> \> a compact open subgroup of $G$ \\
$H^{\perp}$ \> \> the annihilator of $H$ in $\widehat{G}$ \\
$\mu$ \> \> Haar measure on $G$, normalized so that $\mu(H)=1$ \\
$\nu$ \> \> Haar measure on $\widehat{G}$,
	normalized so that $\nu(H^{\perp})=1$
\end{tabbing}
Given $x\in G$ and $\gamma\in \widehat{G}$, the pairing
$(x,\gamma)\in \TT=\{z\in\CC : |z|=1\}$
shall denote the action of $\gamma$ on $x$.
$H^{\perp}$, which by definition consists of all $\gamma\in\widehat{G}$
such that $(x,\gamma)=1$ for each $x\in H$, is a compact open
subgroup of $\widehat{G}$.
The quotient groups $G/H$ and $\widehat{G}/H^{\perp}$ are
discrete, with counting measure induced by $\mu$ and $\nu$,
respectively.  Furthermore, the discrete group
$\widehat{H}$ is naturally
isomorphic to $\widehat{G}/H^{\perp}$, and
the compact group $\widehat{G/H}$ is naturally
isomorphic to $H^{\perp}$.
The normalizations of $\mu$ and $\nu$ are compatible
in the sense that the Fourier inversion formula holds;
see \cite[\S 31.1(c)]{HR2}.

The study of wavelets in $L^2(\RR^d)$ requires a dilation operator
and a discrete group of translation operators.  For the
translations, one usually considers a
discrete lattice $\Lambda\subset\RR^d$ with a fundamental domain
$F$ which is also a neighborhood of the origin.  However, our
group $G$ may not actually contain a discrete subgroup.
Fortunately, the compact subgroup $H$ can play the role
of the fundamental domain, even without an actual lattice.
Meanwhile, the discrete group $G/H$ shall act as a group of
translation operators, in a manner to be described
in Section~\ref{sect:oper}.

The dilation over $\RR^d$ is usually an expansive matrix $B$;
as a linear map, $B$ satisfies $d(\lambda\circ B) = |\det B | d\lambda$,
where $\lambda$ denotes Lebesgue measure on $\RR^d$.
In our setting,
let $A$ be an automorphism of $G$; that is, $A$ is an algebraic isomorphism
from $G$ to $G$ such that $A$ is also a topological homeomorphism.
Then $\mu\circ A$ is a nontrivial Haar measure on $G$.
Define the {\em modulus} of $A$ to be
$|A|=\mu(AH)>0$, so that $d\mu(Ax)=|A|d\mu(x)$.  It is easy to
verify that $|A^{-1}|=|A|^{-1}$.
Furthermore, the adjoint $A^*$ of $A$ is an automorphism of $\widehat{G}$,
and $|A^*|=|A|$ in the sense that $\nu(A^* H^{\perp})=|A|$
and $d\nu(A^*\gamma) = |A| d\nu(\gamma)$.

If the Fourier transform $\widehat{f}$ of
a function $f\in L^2(G)$ has
compact support, then as usual, $f$ is smooth.
Proposition~\ref{prop:locconst} below shall make
that notion of smoothness more precise.  As a preface
to that result, note that
for any fixed $r\in\ZZ$, $A^{-r}H$ is an open subgroup of $G$, and
therefore the sets of the form $c+A^{-r} H$ form a partition
of $G$ into open sets.
The proposition shall show that if the support of $\widehat{f}$
is contained in the compact set $(A^*)^r H^{\perp}$,
then $f$ is constant on all the sets $c+A^{-r} H$.  Thus,
$f$ is locally constant.

\begin{prop}
\label{prop:locconst}
Let $G$ be a locally compact abelian group with compact
open subgroup $H$, let $A$ be an
automorphism of $G$,
let $r\in \ZZ$,
and let $f\in L^2(G)$.

Then  $\Supp \widehat{f} \subset (A^*)^r H^{\perp}$ if and only if
$f$ is constant on every (open) set of the form $c+A^{-r}H$.
\end{prop}

\begin{proof}
To prove the forward implication, pick
any $c\in G$ and $x\in c+ A^{-r}H$.
Note that $x-c\in A^{-r}H$, and therefore
$(x-c,\gamma)=1$ for every $\gamma\in (A^*)^r H^{\perp}$.
Thus,
$$
f(x)
=
\int_{(A^*)^r H^{\perp}}
(x,\gamma) \widehat{f}(\gamma) d\nu(\gamma)
=
\int_{(A^*)^r H^{\perp}}
(c,\gamma) \widehat{f}(\gamma) d\nu(\gamma)
=
f(c).
$$

For the converse, pick $c\in G$ and $r\in\ZZ$.  We shall
show that the desired implication holds for the characteristic
function $f=\mathbf{1}_{c+A^{-r}H}$; the general result shall
then follow by linearity.  We compute
\begin{align}
\label{eq:charhat}
\widehat{f}(\gamma)
&=
\int_{c + A^{-r}H} \bar{(x,\gamma)} d\mu(x)
=
\bar{(c,\gamma)}\int_{A^{-r}H} \bar{(x,\gamma)} d\mu(x)
\notag \\
&=
\bar{(c,\gamma)}|A|^{-r}\int_{H} \bar{(A^{-r}x,\gamma)} d\mu(x)
\notag \\
&=
\bar{(c,\gamma)}|A|^{-r}\int_{H} \bar{(x,(A^*)^{-r}\gamma)} d\mu(x)
\notag \\
&=
\bar{(c,\gamma)}|A|^{-r}\mathbf{1}_{H^\perp}((A^*)^{-r}\gamma)
=
\bar{(c,\gamma)}|A|^{-r}\mathbf{1}_{(A^*)^rH^\perp}(\gamma),
\end{align}
which has support contained in the required domain.
\end{proof}

\begin{definition}
\label{def:expan}
Let $G$ be a locally compact abelian group with compact open
subgroup $H$.  Let $A$ be an automorphism of $G$.  We say
that $A$ is {\em expansive} with respect to $H$ if both
of the following conditions hold:
\begin{enumerate}
\item $AH \supsetneq H$, and
\item $\bigcap_{n\leq 0} A^n H = \{0\}$.
\end{enumerate}
\end{definition}
The idea of Definition~\ref{def:expan} is that an expansive
automorphism acts on $G$ like an expansive integer matrix $B$
acts on $\RR^d$;
the fundamental domain $F$ for the lattice $\ZZ^d$ is contained in
its image $BF$, and the intersection of all inverse images $B^{-k}F$
is just the origin.
Moreover, if $A$ is an expansive automorphism of $G$,
then $|A|$ is an integer greater than $1$, just as is true of an
expansive integer matrix.
Indeed, $H$ is a proper subgroup of $AH$, so that $AH$ may be
covered by (disjoint) cosets $s+H$.  Each coset has measure $1$,
so that $AH$ must have measure equal to the number of cosets.
Thus, $|A|=\mu(AH)$ is just the
number of elements in the finite quotient group $(AH)/H$.

The following alternate characterization of expansiveness
was proven in \cite{BB}.
\begin{prop}[From \cite{BB}]
\label{prop:expan}
Let $G$ be a locally compact abelian group with compact open
subgroup $H$.  Let $A$ be an automorphism of $G$.  Then
$A$ is expansive with respect to $H$ if and only if
both of the following conditions hold:
\begin{enumerate}
\item $A^*H^{\perp} \supsetneq H^{\perp}$, and
\item $\bigcup_{n\geq 0} (A^*)^n H^{\perp} = \widehat{G}$.
\end{enumerate}
\end{prop}

\begin{example}
\label{ex:Qp}
Let $p\geq 2$ be a prime number.
Let $G=\Qp$ be the set of $p$-adic rational numbers,
and $H=\Zp$ the set of $p$-adic integers.
That is,
\begin{align*}
G &= \Qp = \left\{\sum_{n\geq n_0}c_n p^n : n_0\in\ZZ,
c_n=0,1,\ldots, p-1 \right\}, \quad\mbox{and} \\
H &= \Zp = \left\{\sum_{n\geq 0}c_n p^n : c_n=0,1,\ldots, p-1 \right\}.
\end{align*}
Here, the sums are formal sums, but they indicate that the group
law is addition with carrying of digits, so that in $\QQ_3$, for example,
we have
\begin{multline*}
(2 + 2\cdot 3 + 2\cdot3^2 + 2\cdot 3^3 + 2\cdot 3^4 + \ldots)
+ (1 + 2\cdot 3 + 1\cdot3^2 + 2\cdot 3^3 + 1\cdot 3^4 + \ldots)
\\
= (0 + 2\cdot 3 + 1\cdot3^2 + 2\cdot 3^3 + 1\cdot 3^4 + \ldots).
\end{multline*}
The discrete quotient $G/H = \Qp/\Zp$ is naturally isomorphic to
$\mu_{p^{\infty}}$, the multiplicative
group of all $p^n$-th roots of unity
(as $n$ ranges through all nonnegative integers) in $\CC$.
However, $\Qp$ itself contains no nontrivial discrete subgroups;
if a closed subgroup $E$ contains some nonzero element $x\in\Qp$,
then $E$ contains the open subgroup $x\Zp$.

$\Qp$ is in fact a field, with multiplication given by multiplication
of formal Laurent series together with carrying of digits, though we
are considering it only as an additive group.  For more of the standard
background on $\Qp$, see \cite{Gou,Kob,RV}; see \cite{Rob,Ser} for more
advanced expositions.

Define a character $\chi$ on $G$ by
$$\chi\left(\sum_{n\geq n_0}c_n p^n\right) =
\exp\left(2\pi i \sum_{n=n_0}^{-1} c_n p^n\right).$$
Note that $H=\{x\in G: \chi(x)=1\}$.
$G$ is self-dual, with pairing given by
$(x,\gamma) = \chi(x\gamma)$, where $x\gamma$ means the product
in the field $\Qp$.  Under this pairing, $H^{\perp}\subset \Qp$
is again just $\Zp$.

For any nonzero $a\in\Qp$, the multiplication-by-$a$ map
$A:x\mapsto ax$
is an automorphism of $G$.  (In fact, any automorphism must be
of this form.)  If $a\not\in\Zp$, then $A$ is expansive; in
that case, if $a=c_{-m} p^{-m} + c_{-m+1} p^{-m + 1} + \ldots$
with $c_{-m}\neq 0$,
then $|A| = p^m$.
(For example, multiplication-by-$1/p$ is
an expansive map with modulus $p$.)
The adjoint $A^*$ on $\widehat{G}=\Qp$
is again multiplication-by-$a$.
\end{example}

\begin{example}
\label{ex:FpT}
Let $p\geq 2$ be a prime number.
Let $\FF_p$ denote the field of order $p$, with elements
$\{0,1,\ldots, p-1\}$.
Let $G=\FF_p((t))$, the set of formal Laurent series in the variable $t$
with coefficients in $\FF_p$, and let $H=\FF_p[[t]]$ be the subset
consisting of power series (i.e., nonnegative powers of $t$ only).
$G$ is a field under addition and multiplication of formal power series.
Although $G$ is topologically homeomorphic to $\Qp$, the algebraic
structure is very different because there is no carrying of digits.
In $\FF_3((t))$, for example, we have
\begin{multline*}
(2 + 2\cdot t + 2\cdot t^2 + 2\cdot t^3 + 2\cdot t^4 + \ldots)
+ (1 + 2\cdot t + 1\cdot t^2 + 2\cdot t^3 + 1\cdot t^4 + \ldots)
\\
= (0 + 1\cdot t + 0\cdot t^2 + 1\cdot t^3 + 0\cdot t^4 + \ldots).
\end{multline*}
Again, we consider $G$ and $H$ as groups under addition; note
that for any $x\in G$,
$$\underbrace{x + x + \ldots + x}_{p} = 0.$$
(That is, $\FF_p((t))$ is a field of characteristic $p$.)
The discrete quotient $G/H$ is naturally isomorphic to a
direct sum of a countable number of copies of $\FF_p$.
In fact, $G/H$ is isomorphic to the subgroup $\Gamma\subset G$
consisting of all elements of the form
$c_{-m}t^{-m} + \cdots + c_{-1} t^{-1}$.
See \cite{RV} or \cite{Ser} for more on such fields.

Define a character $\chi$ on $G$ by
$$\chi\left(\sum_{n\geq n_0}c_n t^n\right) =
\exp\left(2\pi i c_{-1}/p\right).$$
Note that the image of $\chi$ consists only of the $p$-th roots
of unity, whereas the corresponding character on $\Qp$ included
all $p^n$-th roots of unity in its image.
Another contrast with $\Qp$ is that this time,
$H\subsetneq\{x\in G: \chi(x)=1\}$.
Still, as before, $G$ is self-dual, with pairing given by
$(x,\gamma) = \chi(x\gamma)$, where $x\gamma$ means the product
in the field $\FF_p((t))$.  Under this pairing,
$H^{\perp}\subset \FF_p((t))$
is again just $\FF_p[[t]]$.

For any nonzero $a\in\FF_p((t))$, the
multiplication-by-$a$ map $A:x\mapsto ax$
is an automorphism of $G$.  (This time, however, many other
automorphisms are possible.)
If $a\not\in\FF_p[[t]]$, then $A$ is expansive; in
that case, if $a=c_{-m} t^{-m} + c_{-m+1} t^{-m + 1} + \ldots$
with $c_{-m}\neq 0$,
then $|A| = p^m$.
(For example, multiplication-by-$1/t$ is
an expansive map with modulus $p$.)
The adjoint $A^*$ on $\widehat{G}=\FF_p((t))$
is again multiplication-by-$a$.
\end{example}

\begin{example}
\label{ex:Qpextn}
Let $K$ be a field which is a finite extension of
$\Qp$.
In general, $\Qp$ is equipped with an absolute value
function $|\cdot|:\Qp \rightarrow \RR_{\geq 0}$ with
$|p| = 1/p$ and
which extends uniquely to the algebraic closure of $\Qp$,
and hence to $K$.
The ring of integers $\ints_K$ is precisely the set
of elements of $K$ with absolute value at most $1$,
just as $\ZZ_p$ is the set of elements in $\Qp$
with absolute value at most $1$.
Equivalently, $\ints_K$
consists of all elements of $K$ which are roots
of monic polynomials with coefficients in $\ZZ_p$.

For example, if
$K=\QQ_3(\sqrt{-1})$, then elements of $K$ may
be represented as formal Laurent series $\sum c_n 3^n$,
where each $c_n$ is one of the nine numbers $a+b\sqrt{-1}$,
with $a,b=0,1,2$.  (The addition law for carrying digits
requires some extra computation in this case.)
Similarly, if
$K=\QQ_3(\sqrt{3})$, then elements of $K$ may
be represented as formal Laurent series $\sum c_n 3^{n/2}$,
where each $c_n$ is one of the three numbers $0,1,2$,
and with carrying from the $3^{n/2}$ term to the
$3^{(n+2)/2}$ term.
In both of these cases, $\ints_K$ consists of those
elements of $K$ with no negative power terms.

Let $G$ be the additive group of $K$, and $H\subset G$
be the additive group of $\ints_K$.
Define a character $\chi_K$ on $G$ by
$$\chi_K(x) = \chi(\Tr x),$$
where $\chi$ was the character defined in Example~\ref{ex:Qp}
for $\Qp$, and $\Tr$ is the trace map from $K$ down to $\Qp$.
$G$ is self-dual, with pairing given by
$(x,\gamma) = \chi_K(x\gamma)$, where $x\gamma$ means the product in $K$.
This time, $H^{\perp}\subset K$ need not be the same as $H$.
Instead, $H^{\perp}$ is the inverse different of $K$, which is
a (possibly larger) subgroup of $K$
of the form $b \ints_K$, for some
$b \in K$.  (See \cite[\S B.2]{RV}
or \cite[III]{Ser} for more on the different and inverse
different.)  For example, if $K=\QQ_3(\sqrt{-1})$, then
the inverse different is just $\ints_K =\ZZ_3[\sqrt{-1}]$;
but if $K=\QQ_3(\sqrt{3})$, then the inverse different
is $3^{-1/2} \ints_K = 3^{-1/2} \ZZ_3[\sqrt{3}]$.
(The inverse different roughly measures the
degree of {\em ramification} in $K$, which, speaking informally,
means the largest denominator that appears in the powers of $p$
in expansions of elements of $K$.  So $\QQ_3(\sqrt{3})$ has
ramification degree $2$, since its elements involve half-integral
powers of $3$.)

As before, for any nonzero $a\in K$,
the multiplication $A:x\mapsto ax$
is an automorphism of $G$.  (Other automorphisms are possible,
incorporating Galois automorphisms of $K$ over $\Qp$.)
If $a\not\in\ints_K$, then $A$ is expansive; the
modulus $|A|$ is, as usual, the order of the group $(AH)/H$,
which is always be a power of $p$.

One could also choose $G$ to be the additive group of a field
which is a finite extension of $\FF_p((t))$.  Again, $G$ would
be self-dual, with pairing given by composing the original
$\chi$ for $\FF_p((t))$ with the trace map, and with $H^{\perp}$
given by the inverse different.  We omit the
details here.
\end{example}

\begin{example}
\label{ex:products}
Let $G_1,\ldots, G_N$ be locally compact abelian groups
with compact open subgroups $H_1,\ldots, H_N$,
and expansive automorphisms $A_1,\ldots, A_N$,
respectively.
Then $G=G_1\times\cdots\times G_N$ has open compact subgroup
$H=H_1\times\cdots\times H_N$ with expansive automorphism
$A=A_1\times\cdots\times A_N$.

As a special case, if $G_1=\cdots=G_N=K$ is the additive group of
one of the fields in the previous three examples, then
$G=K^{\times N}$ is a locally compact abelian group with
compact open subgroup $H=\ints_K^{\times N}$.  Moreover, if
$A\in GL(N,K)$ is a matrix with all eigenvalues of
absolute value greater than $1$ and such that $AH\subsetneq H$,
then $A$ is expansive.
The modulus of $A$ is $|\det A |$, where $|\cdot|$ is a certain
appropriately chosen absolute value on $K$.
\end{example}

\begin{example}
\label{ex:weird}
The preceding examples are of self-dual $\sigma$-compact
groups; the following example is of a group which is not
self-dual and may not be $\sigma$-compact.

Let $G_0$ be a locally compact abelian group with
compact open subgroup $H_0$ and expansive automorphism $A_0$,
and let $W$ be any discrete
group (possibly uncountable).  Let $G=G_0\times W$,
with compact open subgroup $H=H_0\times\{0\}$.  Then
$\widehat{G}$ is isomorphic to $\widehat{G}_0 \times \widehat{W}$,
and the image of $H^{\perp}$ under that isomorphism
is $H_0^{\perp} \times \widehat{W}$.
If $W$ is uncountable, note that $G$ is not
$\sigma$-compact.

It is easy to check that the automorphism
$A$ of $G$ given by $A=A_0\times\mbox{id}_W$ (where $\mbox{id}_W$
is the identity automorphism of $W$) is expansive, with
adjoint $A^*=A_0^*\times\mbox{id}_{\widehat{W}}$.  However,
it should be noted that while $\bigcap_{n\leq 0} A^n H = \{0\}$,
as required for expansiveness, the union
$\bigcup_{n\geq 0} A^n H$ is {\em not} all of $G$, but only
$G_0\times \{0\}$.  Similarly, on the dual side,
the union
$\bigcup_{n\geq 0} (A^*)^n H^{\perp}=\widehat{G}$
is the whole group,
but the intersection
$\bigcap_{n\leq 0} (A^*)^n H^{\perp}= \{0\}\times\widehat{W}\neq \{0\}$
includes more than just the identity.
\end{example}

\section{Dilation and translation operators}
\label{sect:oper}
We are now prepared to present the operators on $L^2(G)$
introduced in \cite{BB}.
Given our group $G$ and automorphism $A$,
we can define a
dilation operator $\delta_A$ on $L^2(G)$ by
$$\delta_A(f)(x) = |A|^{1/2} f(Ax).$$
It is easy to verify
that $\delta_A$ is a unitary operator.
If $A$ is expansive, then certainly $A$ does not have finite
order, and therefore $\{\delta_A^n : n\in\ZZ\}$ forms
a group isomorphic to $\ZZ$.

The construction of translation operators requires a little more
work.  If $G$ contains a discrete subgroup $\Gamma$ such that
$G\cong\Gamma\times H$, then $\Gamma$ plays the same role
that the lattice $\ZZ^d$ plays in $\RR^d$.  That is, for any
$s\in\Gamma$, we define $\tau_s$ on $L^2(G)$ by
$\tau_s(f)(x) = f(x-s)$.  For example, if $G=\Fp((t))$,
then the choice of $\Gamma$ described in Example~\ref{ex:FpT}
is such a lattice.  Indeed, Lang \cite{Lan} used precisely this
lattice to define wavelets on $\FF_2((t))$, the Cantor dyadic group.

On the other hand, if $G$ contains no such lattice, as is the
case for $\Qp$, the situation is a little more complicated.
Kozyrev's \cite{Koz} solution for the special case of $\Qp$
was to let $\calC\subset\Qp$ be the discrete subset (but not subgroup)
of elements of the form $c_{-m}p^{-m} + \cdots + c_{-1} p^{-1}$.
Then, for each $s\in\calC$, the translation operator
$\tau_{s}$ is defined in the usual way.  He was then able
to construct Haar wavelets on $\Qp$ using dilation-by-$1/p$
and translation by elements of $\calC$.

More generally, one could choose a set $\calC$ of
coset representatives for the quotient $G/H$.  That is,
choose a subset $\calC\subset G$ such that the translations
$s+H$ of $H$ by elements $s\in\calC$ are mutually disjoint,
and such that $\bigcup_{s\in\calC} (s + H) = G$.  (Such a
$\calC$ always exists if one accepts the Axiom of Choice;
but in many cases, it is possible to construct such a set $\calC$
without invoking the Axiom of Choice.  See, for example,
Proposition~\ref{prop:Ddef}.)

Unfortunately, Kozyrev's method has its limitations.  For
one thing, the resulting group of translation operators
does not form a group unless $\calC$ itself is a group,
which would be impossible for $\Qp$.  Furthermore,
it appears to be difficult or impossible
to construct wavelets other than Haar wavelets
using such translation operators.
More precisely, it appears that
any sort of analogue of a multiresolution analysis or
of minimally supported frequency wavelets requires
the translation operators to behave well on the dual group
$\widehat{G}$, as follows.
If $\calC=\Gamma$ is indeed a group, then
the dual lattice $\Gamma^{\perp}\subset \widehat{G}$
plays a crucial role in wavelet analysis; but if $\calC$
is not a group, then there does not seem to be any reasonable
object which can assume the role of the dual lattice.

In \cite{BB}, the authors solved this problem by
choosing a set $\calD\subset\widehat{G}$ to play the role
of the dual lattice first, and then building the
operators to respect $\calD$.  Specifically, let $\calD$ be a choice
of coset representatives for the quotient $\widehat{G}/H^{\perp}$.
Define functions
$\theta=\theta_{\calD}: \widehat{G} \rightarrow \calD$
and
$\eta=\eta_{\calD}: \widehat{G} \rightarrow H^{\perp}$
given by
\begin{align}
\theta(\gamma) &= \mbox{ the unique } {\sigma} \in\calD
        \mbox{ such that }
	\gamma-\sigma \in H^{\perp}, \notag \\
\eta(\gamma) &= \gamma - \theta(\gamma).
\label{eq:etadef}
\end{align}
Speaking informally, $\eta(\gamma)$ is the unique point in the
fundamental domain $H^{\perp}$ which is congruent to $\gamma$
modulo the ``lattice'' $\calD$, and
$\theta(\gamma)$ is the unique ``lattice'' point such 
that $\gamma = \eta(\gamma) + \theta(\gamma)$.

Any element $[s]\in G/H$ of the discrete quotient is, by definition,
a coset of the form $s+H$, with $s\in G$.
We wish to define a translation operator
associated with $[s]$ which depends only on the coset, not the particular
representative $s$.  Note that on the transform side, the
translation-by-$s$ operator sends $\widehat{f}(\gamma)$ to
$\bar{(s,\gamma)}\widehat{f}(\gamma)$.  Thus, we define
the translation-by-the-coset operator $\tau_{[s]}$ by defining
its action $\tilde{\tau}_{[s]}$ on $L^2(\widehat{G})$ as follows:
\begin{equation}
\label{eq:taudef}
\tilde{\tau}_{[s]}(\widehat{f})(\gamma)
= \bar{(s,\eta(\gamma))}\widehat{f}(\gamma).
\end{equation}
It is easy to check that $\tilde{\tau}_{[s]}$,
and hence $\tau_{[s]}$, are unitary operators.

To simplify notation, define
$$w_{[s]}(\gamma)=\bar{(s,\eta(\gamma))},$$
which is a function in $L^{\infty}(\widehat{G})$, since
$|w_{[s]}(\gamma)|=1$ for all $\gamma\in\widehat{G}$.
We compute
\begin{equation}
\label{eq:wgroup}
w_{[s+t]}(\gamma) = \overline{(s + t , \eta(\gamma))}
= \overline{(s, \eta(\gamma))}\overline{(t , \eta(\gamma))}
= w_{[s]}(\gamma) w_{[t]}(\gamma).
\end{equation}
If we consider $t\in H$, so that $w_{[t]}$ is identically $1$,
then equation~\eqref{eq:wgroup} shows that $w_{[s+t]}=w_{[s]}$.
Thus, as promised,
$w_{[s]}$ and hence $\tau_{[s]}$ depend
only on $[s]$, not on the particular representative $s$.
Equation~\eqref{eq:wgroup} also shows that
the set of operators
$$\Gamma_{G/H} = \{\tau_{[s]} : [s]\in G/H \}$$
forms a group isomorphic to $G/H$.

On $L^2(G)$, the action of $\tau_{[s]}$
is convolution by the pseudo-measure
$\check{w}_{[s]}$, defined to be the inverse transform of
$w_{[s]}$.
By way of comparison,
the usual translation-by-$s$ operator is convolution by a delta measure
centered at $s$.  Proposition~\ref{prop:tauact} below computes
the action of our operator $\tau_{[s]}$, and the result shows
that our pseudo-measure $\check{w}_{[s]}$ has support
contained in $s+H$,
in the sense that it is trivial
on any measurable set
$U\subset G$ disjoint from $s+H$.  Thus, $\check{w}_{[s]}$
may be thought of as a perturbation of the original delta measure.
To see this, we consider the action of $\tau_{[s]}$ on
the characteristic function of a set of the form
$c+A^{-r} H$ in the following proposition.

\begin{prop}
\label{prop:tauact}
Let $G$ be a locally compact abelian group with compact
open subgroup $H$, and let $\calD$ be a set of coset representatives
for the quotient $\widehat{G}/H^{\perp}$.  Let $A$ be an expansive
automorphism of $G$, and let $M=|A|\geq 2$ be the modulus of $A$.
For any $[s]\in G/H$,
define the operator $\tau_{[s]}$ as in equation~\eqref{eq:taudef}
above.  For any $c\in G$ and any $r\in\ZZ$, let
$\mathbf{1}_{c + A^{-r} H}$
be the characteristic function
of the set $c+A^{-r} H$.  
If $r\leq 0$, let $\sigma_0$ be the unique element of $\calD$ in
$H^{\perp}$.
If $r > 0$, let $\{\sigma_0,\sigma_1,\ldots, \sigma_{M^{r}-1}\}$
be all the elements of $\calD\cap (A^*)^{r} H^{\perp}$.
\begin{enumerate}
\item If $r\leq 0$, then
	$$\tau_{[s]}\mathbf{1}_{c + A^{-r} H}(x) =
	(s,\sigma_0) \mathbf{1}_{s + c + A^{-r} H}(x).$$
\item If $r> 0$, then
	$$ \tau_{[s]}\mathbf{1}_{c + A^{-r} H}(x) =
	\frac{1}{M^{r}}
	\left( \sum_{i=0}^{M^{r}-1} (x-c,\sigma_i) \right)
	\mathbf{1}_{s + c + H}(x).$$
\end{enumerate}
\end{prop}

Note that in the case $r\leq 0$ (in which case $A^{-r} H$ is large),
$\tau_{[s]}$ acts as translation-by-$s$ and
multiplication by the constant $(s,\sigma_0)$ on 
$\mathbf{1}_{c + A^{-r} H}$.  For $r> 0$ (in which case
$A^{-r} H$ is small), $\tau_{[s]}$ is similar to translation by $s$,
but it spreads the function out over the full coset $s+c+H$,
not the smaller $s+c+A^{-r} H$.  Still, the translated function
is constant on any region of the form $b+A^{-r} H$.  Moreover,
when that constant is nonzero (i.e., when $b\in s+c+H$), its
value $M^{-r} \sum_i (b,\sigma_i)$ usually is small, due
to a lot of cancellation (but usually not complete cancellation)
within the sum; see the examples later in this section.
Note also that we can easily verify directly from the
proposition that $\tau_{[s]}$ is independent of the coset
representative $s$, and that it is the identity operator
if $s\in H$.

If a function $f\in L^2(G)$ is constant on all such sets
$c+A^{-r} H$, then $f$ 
is a linear combination of characteristic functions
$\mathbf{1}_{c+A^{-r}H}$.
Thus, $\tau_{[s]} f$ may be computed
using Proposition~\ref{prop:tauact} and linearity.
In particular,
the action of $\tau_{[s]}$
on such a function $f$ is determined
not by all of $\calD$ but only by a subset of size $M^r$.
Furthermore, the computation of $\tau_{[s]}f$ involves
only a finite sum and is therefore easily computable in practice.

\begin{proof}[Proof of Proposition~\ref{prop:tauact}]
Let $f=\mathbf{1}_{c + A^{-r}H}$.
By Equation~\eqref{eq:charhat},
$$\widehat{f}(\gamma) = M^{-r} \bar{(c,\gamma)}
	\mathbf{1}_{(A^*)^{r} H^{\perp}} (\gamma),$$
and therefore
$$\tilde{\tau}_{[s]} \widehat{f}(\gamma)
	= M^{-r} \bar{(s,\eta(\gamma))}
	\bar{(c,\gamma)} \mathbf{1}_{(A^*)^{r} H^{\perp}} (\gamma).$$
If $r\leq 0$, then $(A^*)^{r} H^{\perp}\subset H^{\perp}$,
so that $\eta(\gamma)=\gamma-\sigma_0$ for
$\gamma\in (A^*)^{r} H^{\perp}$.  In that case,
$$\tilde{\tau}_{[s]} \widehat{f}(\gamma)
	= M^{-r} (s,\sigma_0)
	\bar{(s+c,\gamma)} \mathbf{1}_{(A^*)^{r} H^{\perp}} (\gamma).$$
Taking the inverse transform gives
$$\tau_{[s]}f(x) = M^{-r} (s,\sigma_0) \int_{(A^*)^{r}}
	(x-s-c,\gamma) d\nu(\gamma)
	= (s,\sigma_0) \mathbf{1}_{s+c+ A^{-r} H}(x).$$

On the other hand, if $r>0$, then
$(A^*)^{r} H^{\perp}\supsetneq H^{\perp}$.  To compute
$\eta(\gamma)$, then, we partition
$(A^*)^{r} H^{\perp}$ into the
$M^{r}$ regions $\{\sigma_i + H^{\perp} : i=0,1,\ldots, M^{r}-1\}$.
We have $\eta(\gamma)=\gamma-\sigma_i$ for
$\gamma\in\sigma_i + H^{\perp}$.  Thus,
$$\tilde{\tau}_{[s]} \widehat{f}(\gamma)
	= M^{-r} \sum_{i=0}^{M^{r}-1}
	(s,\sigma_i) \bar{(s+c, \gamma)}
	\mathbf{1}_{\sigma_i +  H^{\perp}} (\gamma).$$
Taking the inverse transform,
\begin{align*}
\tau_{[s]}f(x) &= M^{-r} \sum_{i=0}^{M^{r}-1}
(s,\sigma_i) \int_{\sigma_i + H^{\perp}}
(x-s-c,\gamma) d\nu(\gamma)
\\
&= M^{-r} \sum_{i=0}^{M^{r}-1}
(x-c,\sigma_i) \mathbf{1}_{s+c+H}(x). \qedhere
\end{align*}
\end{proof}

Of course, the operator $\tau_{[s]}$ depends on the choice $\calD$ of
coset representatives.  In \cite{BB} it was shown that
given an expansive automorphism $A$ of $G$,
it is possible to construct such a $\calD$ while making only
finitely many choices.  We restate that result here as
a proposition.
\begin{prop}[From \cite{BB}]
\label{prop:Ddef}
Let $G$ be a locally compact abelian group with compact
open subgroup $H$.  Let $A$ be an expansive automorphism of $G$
with modulus $|A|=M\geq 2$,
and let $\calD_1 = \{\rho_0,\ldots,\rho_{M-1}\}$ be
a set of coset representatives for the (finite) quotient
$H^{\perp}/((A^*)^{-1}H^{\perp})$, with $\rho_0=0$.
Define $\calD \subset \widehat{G}$ to be the
(infinite) set of
all elements $\sigma\in\widehat{G}$ of the form
\begin{equation}
\label{eq:Ddef}
\sigma = \sum_{j=1}^n (A^*)^j \rho_{i_j},
\quad \text{where } n\geq 1 \text{ and } i_j\in\{0,1,\ldots, M-1\}.
\end{equation}
Then $\calD$ is a set of coset representatives for the
quotient $\widehat{G} / H^{\perp}$.
Moreover, $A^* \calD\subsetneq \calD$.
\end{prop}

\begin{example}
\label{ex:QpD}
Let $p\geq 2$ be a prime number, $G=\Qp$, and $H=\Zp$;
identify $\widehat{G}$ as $\Qp$ and $H^{\perp}$ as $\Zp$,
as in Example~\ref{ex:Qp}.
Let the automorphism $A$ be multiplication-by-$1/p$.
As mentioned previously, Kozyrev \cite{Koz} chose
a certain standard set of coset representatives in $\Qp$ for
the quotient $G/H$; the same set in $\Qp$ can be used
as our choice of coset representatives for $\widehat{G}/H^{\perp}$.
Specifically, let $\calD$ consist of all elements of
$\widehat{G}=\Qp$ of the form
$$\sigma = \sum_{i=1}^n c_i p^{-i}, \quad
\text{with } c_i\in\{0,1,\ldots, p-1\}.$$
This is the same set $\calD$ that would come from
$\calD_1=\{0,1,2\ldots, (p-1)\}$ using
Proposition~\ref{prop:Ddef}.

The resulting operator $\tau_{[s]}$ takes $f=\mathbf{1}_{c+p^n\Zp}$
to
$$\tau_{[s]}f(x)=
p^{-n}\left[\sum_{j=0}^{p^{n}-1} (x-c,\sigma_j)\right]
\mathbf{1}_{s+c+\Zp}(x)$$
by Proposition~\ref{prop:tauact}, with $r=n>0$.
If $s\in\Zp$ (that is, if $[s]=[0]$), then it is easy to
verify from the above formula that $\tau_{[s]}f=f$.
On the other hand, if $s\not\in\Zp$,
write $x\in s+c+\Zp$ as
$x\in c+ m/p^{\ell} + p^{n}\ZZ_p$,
for some integer $m$ prime to $p$ and some $\ell>1$.
(Note that $m$ is only determined modulo $p^{n+\ell}$.)
The sum in the formula for $\tau_{[s]}f$ becomes
$$p^{-n}\left[\sum_{j=0}^{p^{n}-1} (x-c,\sigma_j)\right]
=p^{-n}\sum_{j=0}^{p^{n}-1} (m/p^{\ell},j/p^n)
=p^{-n}\sum_{j=0}^{p^{n}-1} e^{2\pi i mj/p^{\ell+n}}.
$$
We may assume that $1\leq m\leq p^{\ell+n}-1$.  For values
of $m$ larger than $p^{\ell}$, we should expect a great deal
of cancellation in the sum, because the terms
must be distributed fairly evenly around the unit circle.
Thus, most of the weight of $\tau_{[s]}\mathbf{1}_{c+p^n\Zp}$
should be concentrated in the (relatively few) regions $b+p^n\Zp$
corresponding to small values of $m$.

We note that another standard choice for $\calD$ would be to use
$\calD_1=\{0\}\cup\{\zeta_{p-1}^i : i=1,\ldots,(p-1)\}$,
where $\zeta_{p-1}\in\Qp$ is a primitive $(p-1)$st root of unity.
\end{example}

\begin{example}
\label{ex:FpTD}
Let $p\geq 2$ be a prime number, $G=\FF_p((t))$, and $H=\FF_p[[t]]$;
identify $\widehat{G}$ as $\FF_p((t))$ and $H^{\perp}$ as $\FF_p[[t]]$,
as in Example~\ref{ex:FpT}.
Let the automorphism $A$ be multiplication-by-$1/t$.
We may choose our set $\calD$ of coset representatives for
$\widehat{G}/H^{\perp}$ to be all elements of $\widehat{G}$ of the form
$$\sigma = \sum_{i=1}^n c_i t^{-i}, \quad
\text{with } c_i\in\FF_p.$$
This is the same set $\calD$ that would come from
Proposition~\ref{prop:Ddef} using $\calD_1=\FF_p$ (the constants
in $\FF_p[[t]]$).  As a subset of $\FF_p((t))$, it is the
same as the group $\Gamma$
mentioned in Example~\ref{ex:FpT}.  This time we are
considering it as a subgroup of $\widehat{G}$, not of $G$, and
so we call the group $\Gamma'=\calD$.

In particular, $\calD$ is actually a subgroup of $\widehat{G}$.
Thus, the sum
$$p^{-n}\sum_{j=0}^{p^{n}-1} (x-c,\sigma_j)$$
(where $\sigma_j$ ranges over the $p^n$ elements
of the group $\Gamma' \cap t^{-n}\FF_p[[t]]$)
from Proposition~\ref{prop:tauact} simplifies to 
zero unless $(x-c,\sigma_j)$ is trivial for all $j=0,1,\ldots, p^n-1$.
It follows that the operator $\tau_{[s]}$ is the usual
translation-by-$u$ operator, where $u\in G$ is the
unique element of $(s+H)\cap\Gamma$.  In particular,
for $p=2$, in which case $G=\FF_2((t))$ is sometimes known
as the ``Cantor dyadic group'',
the translation operators induced by $\calD$ are the same
as the translation operators used by Lang in his study
of wavelets on that group \cite{Lan}.  Indeed, for any $p$,
$\FF_p((t))$ is a locally compact abelian group with a
discrete lattice $\Gamma$ and dual lattice $\Gamma'$,
which is the situation studied by Dahlke in \cite{Dah}.

We close this example by noting that one could choose $\calD$
in other ways so as not to form a subgroup.  In that case, the
operators $\tau_{[s]}$ would not be true translation operators,
but would behave more like the operators of Example~\ref{ex:QpD}
above.
\end{example}

\begin{example}
\label{ex:QpextnD}
Let $G= \QQ_3(\sqrt{-1})$ and $H=\ZZ_3[\sqrt{-1}]$;
identify $\widehat{G}$ as $\QQ_3(\sqrt{-1})$ and $H^{\perp}$
as $\ZZ_3[\sqrt{-1}]$.
Let $A$ be multiplication-by-$1/3$, so that $|A|=9$.
We may choose $\calD_1$ to be the set of nine elements
$$\calD_1 = \{a+b\sqrt{-1} : a,b =0,1,2\},$$
so that the resulting full set $\calD$ of coset representatives
for $\widehat{G}/H^{\perp}$ consists of all elements of the form
$$ \sigma = \sum_{i=1}^n (a_i + b_i \sqrt{-1}) 3^{-i}, \quad
	\text{with } a_i,b_i\in \{0,1,2\}.$$
Thus, $\calD$ consists of elements of the form $q_1 + q_2\sqrt{-1}$,
with $q_j= m_j/3^{n_j}$ and $0\leq m_j \leq 3^{n_j}-1$.

The resulting operator $\tau_{[s]}$ takes $f=\mathbf{1}_{c+3^n H}$
to
$$\tau_{[s]}f(x)=
9^{-n}\left[\sum_{m_1=0}^{3^{n}-1}\sum_{m_2=0}^{3^{n}-1}
 (x-c,(m_1+m_2\sqrt{-1})/3^n)\right] \mathbf{1}_{s+c+H}(x)$$
by Proposition~\ref{prop:tauact}, with $r=n>0$.
If $s\in H$ (that is, if $[s]=[0]$), then it is easy to
verify from the above formula that $\tau_{[s]}f=f$.
On the other hand, if $s\not\in H$,
write $x\in s+c+\Zp$ as
$x\in c+ (m_3 + m_4\sqrt{-1})/3^{\ell} + 3^{n}H$,
for some $\ell>1$ and some integers $m_3,m_4$ not both divisible
by $3$.
(Note that $m_3$ and $m_4$ are only determined modulo $p^{n+\ell}$.)
The trace map from $\QQ_3(\sqrt{-1})$ to $\QQ_3$ takes
$r_1 + r_2\sqrt{-1}$ to $2r_1$, and so
the sum in the formula for $\tau_{[s]}f$ becomes
\begin{align*}
\tau_{[s]}f(x) & =
9^{-n}\left[\sum_{m_1=0}^{3^{n}-1}\sum_{m_2=0}^{3^{n}-1}
 (x-c,(m_1+m_2\sqrt{-1})/3^n)\right]
\\
{} & = 9^{-n}\sum_{m_1=0}^{3^{n}-1}\sum_{m_2=0}^{3^{n}-1}
e^{4\pi i (m_1 m_3 - m_2 m_4)/3^{\ell+n}}.
\end{align*}
We may assume that $1\leq m_3,m_4\leq p^{\ell+n}-1$.
As in Example~\ref{ex:QpD}, we should expect a great deal
of cancellation in the sum if either or both of $m_3$ and $m_4$
are larger than $p^{\ell}$.
Thus, most of the weight of $\tau_{[s]}\mathbf{1}_{c+3^n H}$
should be concentrated in the (relatively few) regions $b+3^n H$
corresponding to small values of $m_3$ and $m_4$.
\end{example}

\section{Haar and Shannon wavelets}
\label{sect:haar}

Given our dilation operator $\delta_{A}$ and our group of
translation operators $\Gamma_{G/H}=\{\tau_{[s]} \}$, we
may now define wavelets on $G$ as in \cite{BB}.
\begin{definition}
\label{def:wvlt}
Let $G$ be a locally compact abelian group
with compact open subgroup $H \subseteq G.$
Let $A$ be an automorphism of $G$, expansive
with respect to $H$.
Let $\calD$ be a choice of coset representatives in $\widehat{G}$
for $\widehat{H}=\widehat{G}/H^{\perp}$,
and let $\Gamma_{G/H} = \{\tau_{[s]} : [s]\in G/H \}$
be the group of translation operators on $L^2(G)$ determined
by $\calD$ via equation~\eqref{eq:taudef}.

Let $\Psi=\{\psi_1,\ldots,\psi_N\} \subseteq L^2(G)$ be
a finite subset of $L^2(G)$.
$\Psi$ is
a {\em set of wavelet generators} for $L^2(G)$ with respect to $\calD$
and $A$ if
$$
\{\delta_A^n \tau_{[s]} \psi_i : 1\leq i \leq N, n\in\ZZ, [s]\in G/H \}
$$
forms an orthonormal basis for $L^2(G)$.  In that case, the
resulting basis is called a {\em wavelet basis} for $L^2(G)$.

If $N=1$ and $\Psi = \{\psi\}$, 
then $\psi$ is a {\em wavelet} for $L^2(G).$
\end{definition}

The following theorems show that any group $G$ of the type we
have considered has both Haar and Shannon wavelets, and, perhaps
surprisingly, that the Haar wavelets are precisely the same as
the Shannon wavelets.

\begin{theorem}
\label{thm:haardef}
Let $G$ be a locally compact abelian group with compact
open subgroup $H$, and let $\calD$ be a set of coset representatives
for the quotient $\widehat{G}/H^{\perp}$.  Let $A$ be an expansive
automorphism of $G$, and let $M=|A|\geq 2$ be the modulus of $A$.
Let $\{\sigma_0,\sigma_1,\ldots,\sigma_{M-1}\}$ be the $M$ elements
of $(A^*H^{\perp})/H^{\perp}$, with $\sigma_0\in H^{\perp}$.

Define
$$ \Omega_i= H^{\perp}+\sigma_i, \quad \text{and} \quad
\psi_i = \check{\mathbf{1}}_{\Omega_i}, \quad
\text{for all } i=1,\ldots, M-1.$$
That is, $\psi_i$ is the inverse transform of the characteristic
function of $H^{\perp}+\sigma_i$.
Then $\{\psi_1,\ldots,\psi_{M-1}\}$ is a set of wavelet
generators for $L^2(G)$ with respect to $\calD$
and $A$.
\end{theorem}

\begin{theorem}
\label{thm:haarshan}
Let $G$, $H$, $\calD$, $A$, $M$, $\{\sigma_i : i=0,\ldots, M-1\}$,
and $\{\psi_i : i=1,\ldots, M-1\}$ be as in
Theorem~\ref{thm:haardef}.
Then for any $i=1,\ldots, M-1$ and any $[s]\in G/H$,
the translated function $\tau_{[s]}\psi_i$ is
$$\tau_{[s]}\psi_i(x) = (x,\sigma_i) \mathbf{1}_{s+H}(x).$$
Moreover, $\tau_{[s]}\psi_i$ is locally constant;
specifically, it is constant
on every set of the form $c+A^{-1}H$.
\end{theorem}

\begin{proof}[Proof of Theorem~\ref{thm:haardef}]
The sets $(A^*)^n \left(\sigma_i + H^{\perp}\right)$
are pairwise disjoint as $n$ ranges over $\ZZ$ and $i$ ranges
over $\{1,\ldots, M-1\}$.  Moreover, the union of all those
sets is $\widehat{G}\setminus\{0\}$.  Thus, up to sets of
measure zero,
$$\left\{ (A^*)^n \Omega_i : n\in\ZZ, i=1,2,\ldots, M-1 \right\}$$
tiles $\widehat{G}$.
At the same time, each $\Omega_i$ is
the translate of $H^{\perp}$ by $\sigma_i-\sigma_0$.
By the results of \cite{BB}, then, 
$\{\psi_1,\ldots,\psi_{M-1}\}$ is a set of wavelet
generators.
\end{proof}

\begin{proof}[Proof of Theorem~\ref{thm:haarshan}]
By definition, $\tilde{\tau}_{[s]}\widehat{\psi}_i$ is the function
$$\tilde{\tau}_{[s]}\widehat{\psi}_i(\gamma)
= \bar{\left( s, \eta(\gamma) \right)} \widehat{\psi}_i(\gamma)
= \bar{\left( s, \gamma - \sigma_i \right)}
\mathbf{1}_{\sigma_i + H^{\perp}}(\gamma).$$
We compute $\tau_{[s]}\psi_i$ by taking the inverse transform:
\begin{align*}
\tau_{[s]}\psi_i(x)
&=
\int_{\sigma_i + H^{\perp}}
\bar{\left( s, \gamma -\sigma_i \right)}
\left( x, \gamma \right)
d\nu(\gamma)
=
\int_{H^{\perp}}
\bar{\left( s, \gamma \right)}
\left( x, \gamma \right)
\left( x, \sigma_i \right)
d\nu(\gamma)
\\
&=
\left( x, \sigma_i \right)
\int_{H^{\perp}}
\bar{\left( x-s, \gamma \right)}
d\nu(\gamma)
=
\left( x, \sigma_i \right)
\mathbf{1}_{s+H}(x)
\end{align*}
as claimed.

Finally, $\tau_{[s]}\psi_i$ is constant on any set of
the form $c + A^{-1} H$, by Propostion~\ref{prop:locconst}.
This local constancy condition may also be proven directly
from the above formula for $\tau_{[s]}\psi_i$; we leave
the details to the reader.
\end{proof}

By Theorem~\ref{thm:haardef}, the transform
$\widehat{\psi}_i$ of each $\psi_i$ is the characteristic
function of an uncomplicated set $\Omega_i=\sigma_i + H^{\perp}$
formed by translating the fundamental domain $H^{\perp}$.
Moreover, these $M-1$ sets, together with
$\Omega_0 = H^{\perp}$ itself,
tile $A^* H^{\perp}$.  Thus, the functions $\{\psi_i\}$
may be considered analogues for $G$ of the standard
Shannon wavelets from $\RR^d$.

On the other hand, by Theorem~\ref{thm:haarshan},
each $\psi_i$ is a function supported in $H$
and which is constant on
each of the $M$ subsets of $H$ of the form $c+A^{-1}H$.
Furthermore, each of those constant values is an $M$-th root
of unity in $\CC$.  (After all, $M\sigma_i\in H$, since
$AH/H$ is a group of order $M$.)  Therefore, the functions
$\{\psi_i\}$ are analogues for $G$ of the standard Haar
wavelets for $\RR^d$.

As claimed, then,
Shannon wavelets are the same as Haar wavelets
for our group $G$.  Moreover,
these wavelets $\psi_i$ are locally constant and have
compact support; the same is true of their transforms
$\widehat{\psi}_i$.  Thus, our Haar/Shannon wavelets are
simultaneously smooth and of compact
support in both the space and transform domains.

According to Theorem~\ref{thm:haarshan}, the Haar wavelet
$\psi_i$ satisfies 
$$\psi_i(x) = (x,\sigma_i) \mathbf{1}_{H}(x), \quad \text{and} \quad
\tau_{[s]}\psi_i(x) = (x,\sigma_i) \mathbf{1}_{s+H}(x).$$
On the other hand, the usual translation-by-$s$ operator
sends $\psi_i$ to
$$\psi_i(x-s) = \bar{(s,\sigma_i)} (x,\sigma_i) \mathbf{1}_{s+H}(x),$$
which is simply the constant $\bar{(s,\sigma_i)}$ multiplied
by $\tau_{[s]} \psi_i(x)$.  In particular,
as Kozyrev \cite{Koz} showed for the special case of $\Qp$,
$\{\psi_i: 1\leq i\leq M-1\}$
is a set of wavelet generators even if the usual translation-by-$s$
operators are used, as $s$ ranges over a set $\calC$ of coset
representatives for $G/H$.

Similarly, if we replace our ``dual lattice'' $\calD$ by another
choice $\calD'$ of coset representatives for $\widehat{G}/H^{\perp}$,
then the functions $\tau'_{[s]}\psi_i$
are different from $\tau_{[s]}\psi_i$, but only by a
constant multiple $(s,\sigma_i - \sigma'_i)$.
Thus, $\{\psi_i: 1\leq i\leq M-1\}$ is still a set of wavelet
generators for the translation operators $\tau'_{[s]}$ induced
by $\calD'$.  

The fact that waveletness is preserved even when the translation
operators are changed so arbitrarily is unique to these particular
wavelets.  It happens mainly because the support of $\widehat{\psi}_i$
is contained in a single coset $\sigma+H^{\perp}$, and therefore
$\psi_i$ is the exponential $(x,\sigma)$ multiplied by
the characteristic function of $H$.
There are simply not very many candidates for
a new function $\tau_{[s]}\psi_i$ which could reasonably
be considered some kind of translation of $\psi_i$
by the coset $s+H$.

Other wavelets have transforms $\widehat{\psi}$ with
support in at least two cosets $\sigma + H^{\perp}$, and
therefore the usual translation operators do not work as
well.  In particular, the fact that there is no lattice
in $G=\Qp$ (or many other such groups) means that no discrete
set of usual translation operators can form a group,
making it very difficult for all the translations
of a non-Haar $\psi$ to be orthonormal.
In the next section we shall
consider a few examples of such non-Haar wavelets.

\section{Other examples}
\label{sect:exs}

If our expansive map $A$ has modulus $|A|\geq 3$, then
the sets of wavelet generators
described in Section~\ref{sect:haar} have at least
two elements.  Using wavelet sets, J.~Benedetto and
the author produced (single) wavelets in the same
context \cite{BB}.  We state a special case of
that result here without proof:
\begin{theorem}[From \cite{BB}]
\label{thm:bb}
Let $G$ be a locally compact abelian group with compact
open subgroup $H$, let $A$ be an expansive
automorphism of $G$, and let $M=|A|\geq 2$.
Let $\calD$ be a set of coset representatives
for the quotient $\widehat{G}/H^{\perp}$, and
let $\sigma_0$ be the unique element of $\calD \cap H^{\perp}$.
Define
$$\calE = \left(\calD \cap A^* H^{\perp} \right)
	\setminus \{\sigma_0 \},$$
and let $\{V_{\sigma} : \sigma\in \calE\}$
be a partition of $H^{\perp}$ into $M-1$ measurable subsets.

Define a function $T:H^{\perp}\rightarrow \widehat{G}$ by
$$T(\gamma) = \gamma + \sigma -\sigma_0
\quad \text{for } \gamma\in V_{\sigma}.$$
Then there is a measurable set $\Omega\subset A^* H^{\perp}$
such that, up to measure zero subsets,
\begin{equation}
\label{eq:Omegadef}
\Omega =
\left( H^{\perp} \setminus
\bigcup_{n\geq 1} (A^*)^{-n} \Omega \right)
\cup T\left( \bigcup_{n\geq 1} (A^*)^{-n} \Omega \right).
\end{equation}
Let $\psi=\check{\mathbf{1}}_{\Omega}$ be the inverse
transform of the characteristic function of $\Omega$.
Then $\psi$ is a wavelet for $L^2(G)$.
\end{theorem}

The set $\Omega$ may be constructed from $T$ inductively
as follows.  Define $\Omega_0= H^{\perp}$; for every
$n\geq 1$, define
$$\Lambda_n = \Omega_{n-1} \cap \left(
\bigcup_{i\geq 1} (A^*)^{-i} \Omega_{n-1} \right)
\quad\text{and} \quad
\Omega_n = \left( \Omega_{n-1} \setminus \Lambda_n \right)
\cup T\Lambda_n . $$
Then let
\begin{equation}
\label{eq:Omegadef2}
\Lambda = \bigcup_{n\geq 1} \Lambda_n ,
\qquad
\text{and}
\qquad
\Omega = \left( H^{\perp} \setminus \Lambda \right) \cup T\Lambda.
\end{equation}
We may compute the resulting wavelet $\psi$, and any
translate $\tau_{[s]}\psi$, as follows.  With notation
as in Theorem~\ref{thm:bb}, we have
\begin{align*}
\tau_{[s]}\psi(x) & =
\int_{\Omega} (x,\gamma) \bar{(s, \eta(\gamma))} d\nu(\gamma)
\\
& =
\int_{H^{\perp}\setminus \Lambda} (x,\gamma)
\bar{(s, \gamma - \sigma_0)} d\nu(\gamma)
+\sum_{\sigma\in\calE} \int_{T(\Lambda\cap V_{\sigma})}
(x,\gamma) \bar{(s, \gamma - \sigma)} d\nu(\gamma)
\\
& =
(s,\sigma_0) \left[
\int_{H^{\perp}\setminus \Lambda} (x-s,\gamma) d\nu(\gamma)
+\sum_{\sigma\in\calE} \int_{(\Lambda\cap V_{\sigma})}
(x,\gamma+\sigma-\sigma_0) \bar{(s, \gamma)} d\nu(\gamma) \right] .
\end{align*}
If we add and subtract
$(s,\sigma_0)\int_{\Lambda} (x-s,\gamma) d\nu(\gamma)$, we see that
the $\psi$ produced by Theorem~\ref{thm:bb} is given by
\begin{equation}
\label{eq:psicomp}
\tau_{[s]}\psi(x) =
(s,\sigma_0)\left[
\mathbf{1}_H(x-s)
+ \sum_{\sigma\in\calE} \left( (x,\sigma-\sigma_0) - 1 \right)
\int_{(\Lambda\cap V_{\sigma})}
(x-s,\gamma) d\nu(\gamma) \right].
\end{equation}

Note that $\psi$
has transform $\widehat{\psi}=\mathbf{1}_{\Omega}$ with support
contained in $A^* H^{\perp}$.
Thus, for any $[s]\in G/H$, $\tau_{[s]}\psi$ is constant
on any set of the form $c+A^{-1} H$,
by Proposition~\ref{prop:locconst}.

In the special case that $|A|=2$, then $\calE$ has only one element
$\sigma_1$.  Therefore
there is only one choice for $T$ (namely
$T(\gamma)=\gamma+\sigma_1 - \sigma_0$ for all $\gamma\in H^{\perp}$),
which produces $\Omega=\sigma_1 + H^{\perp}$.  The resulting
$\psi$ is simply the Haar/Shannon wavelet from
Section~\ref{sect:haar}.

On the other hand, if $|A|\geq 3$, then although $\psi$ is locally
constant, it does not have compact support.  We now consider
examples of such wavelets $\psi$.

\begin{example}
\label{ex:Qpwave}
Let $p\geq 2$ be a prime number, let $G=\Qp$, let $H=\Zp$, and
let $A$ be multiplication-by-$1/p^r$, for some integer $r\geq 1$;
note that $|A|=p^r$.
Identify $\widehat{G}$ as $\Qp$ and $H^{\perp}$ as $\Zp$,
as in Example~\ref{ex:Qp}.
Let $\calD$ be the set of coset representatives for
$\widehat{G}/H^{\perp}$ defined in Example~\ref{ex:QpD}.
Define $V_{1/p^r}=H^{\perp}$, and $V_{\sigma}=\emptyset$
for all $\sigma\in \calD \setminus \{1/p^r\}$; then define
$T$ as in Theorem~\ref{thm:bb}.  That is, for all $\gamma\in H^{\perp}$,
$T(\gamma) = \gamma + 1/p^r$.  For every $n\geq 1$,
the subset $\Lambda_n$ of $\widehat{G}$ is
$$\Lambda_n = 1 + p^r + p^{2r} + \cdots + p^{(n-2)r} + p^{nr}\Zp,$$
which is the closed ball of radius $1/p^{nr}$ about 
$1+p^r+\cdots+p^{(n-2)r}\in \widehat{G}$,
and which has measure $1/p^{nr}$.
(Note that $\Lambda_1=p^r\Zp$
is the closed ball of radius $1/p^r$ about $0$.)
Define $\Lambda$ and $\Omega$ as in equation~\eqref{eq:Omegadef2}.
Note that $\Lambda=\bigcup \Lambda_n \subset\Zp$
is a disjoint union of countably many
balls, and that $\nu(\Lambda) = 1/(p^r-1)$.
Define $\psi=\check{\mathbf{1}}_{\Omega}$;
then by Theorem~\ref{thm:bb}, $\psi$ is a (single) wavelet.

We can describe $\tau_{[s]}\psi(x)$ explicitly using
equation~\eqref{eq:psicomp}, as follows.  Since there is only
one nonempty $V_{\sigma}$, namely $V_{1/p^r}=\Zp$, we have
$\Lambda\cap V_{1/p^r}=\Lambda$, and therefore
$$\tau_{[s]}\psi(x) =
\mathbf{1}_{\Zp} (x-s) +
\left( (x,1/p^r) - 1\right)
\sum_{n\geq 1}\int_{\Lambda_n} (x-s,\gamma) d\gamma.$$
Given $x,s\in G$ and $n\geq 1$, we compute
$$\int_{\Lambda_n} (x-s,\gamma) d\gamma =
\frac{1}{p^{nr}}(x-s,1+p^r+\cdots+p^{(n-2)r})
\mathbf{1}_{p^{-nr}\ZZ_p}(x-s).$$
Thus, for $x,s\in G$, if we
let $N$ be the smallest positive integer such
that $x-s\in p^{-Nr}\Zp$, then
\begin{align*}
\tau_{[s]}\psi(x) =
\mathbf{1}_{\Zp} (x-s) & +
\frac{\left[(x,1/p^r) - 1\right]}{p^{Nr}} 
(x-s,1+p^r+\cdots+p^{(N-2)r})
\\
& + 
\frac{\left[ (x,1/p^r) - 1\right]}{p^{Nr}(p^r - 1)} 
(x-s,1+p^r+\cdots+p^{(N-1)r}),
\end{align*}
because $(x-s,\gamma)=1$ for all $\gamma\in p^{Nr}\Zp$.
Note that the $p^{Nr}$ in the denominator forces $\tau_{[s]}\psi$
to decay as $x$ goes to $\infty$ in $\Qp$.
\end{example}

\begin{example}
\label{ex:Qpextnwave}
Let $G=\QQ_3(\sqrt{-1})$, let $H=\ZZ_3[\sqrt{-1}]$, and
let $A$ be multiplication-by-$1/3^r$, for some integer $r\geq 1$;
note that $|A|=9^r$.
Identify $\widehat{G}$ as $\QQ_3(\sqrt{-1})$ and
$H^{\perp}$ as $\ZZ_3[\sqrt{-1}]$,
as in Example~\ref{ex:Qpextn}.
Let $\calD$ be the set of coset representatives for
$\widehat{G}/H^{\perp}$ defined in Example~\ref{ex:QpextnD}.
Define $V_{1/3^r}=H^{\perp}$, and $V_{\sigma}=\emptyset$
for all $\sigma\in \calD \setminus \{1/3^r\}$; then define
$T$ as in Theorem~\ref{thm:bb}.  That is, for all $\gamma\in H^{\perp}$,
$T(\gamma) = \gamma + 1/3^r$.  For every $n\geq 1$,
the subset $\Lambda_n$ of $\widehat{G}$ is
$$\Lambda_n = 1 + 3^r + 3^{2r} + \cdots + 3^{(n-2)r}
+ 3^{nr}\ZZ_3[\sqrt{-1}],$$
which is the closed ball of radius $1/3^{nr}$ about 
$1+3^r+\cdots+3^{(n-2)r}\in \widehat{G}$,
and which has measure $1/9^{nr}$.

As before, the set 
$\Lambda$ from equation~\eqref{eq:Omegadef2} is
a disjoint union of countably many
balls, and that $\nu(\Lambda) = 1/(9^r-1)$.
Define $\psi=\check{\mathbf{1}}_{\Omega}$;
then by Theorem~\ref{thm:bb}, $\psi$ is a (single) wavelet.
By a computation similar to that in
Example~\ref{ex:Qpwave}, we find that
\begin{align*}
\tau_{[s]}\psi(x) =
\mathbf{1}_{H} (x-s) & +
\frac{\left[(x,1/3^r) - 1\right]}{9^{Nr}} 
(x-s,1+3^r+\cdots+3^{(N-2)r})
\\
& + 
\frac{\left[ (x,1/3^r) - 1\right]}{9^{Nr}(9^r - 1)} 
(x-s,1+3^r+\cdots+3^{(N-1)r}),
\end{align*}
where $N$ be the smallest positive integer such
that $x-s\in 3^{-Nr} H$.
\end{example}

\begin{example}
\label{ex:FpTwave}
Let $p\geq 2$ be a prime number, let $G=\FF_p((t))$,
let $H=\FF_p[[t]]$, let $A$ be multiplication-by-$1/t$,
and identify $\widehat{G}$ as $\FF_p((t))$ and
$H^{\perp}$ as $\FF_p[[t]]$, as in Example~\ref{ex:FpT}.
Let $\calD$ be the set of coset representatives for
$\widehat{G}/H^{\perp}$ defined in Example~\ref{ex:FpTD}.
Define $V_{1/t}=H^{\perp}$, and $V_{\sigma}=\emptyset$
for all $\sigma\in \calD \setminus \{1/t\}$, and define
$T$ as in Theorem~\ref{thm:bb}.
That is, for all $\gamma\in H^{\perp}$,
$T(\gamma) = \gamma + 1/t$.

The resulting $\Lambda_n$, $\Lambda$, and $\Omega$ are
exactly analogous to those in Example~\ref{ex:Qpwave},
but with $t$ in place of $p$.  As in that case, we get
$$\tau_{[s]}\psi(x) =
\mathbf{1}_{\FF_p[[t]]} (x-s) +
\left[\frac{(x,1/t) - 1}{p^N (p-1)}\right]
(x-s,1+t+\cdots+t^{N-1}).$$
(We still have $p$ appearing in the denominator
because the measure of $\Lambda_n$ is the
{\em real number} $p^n$, not an element of $G$.)

This time, we can simplify the expression even further.
Write
\begin{align*}
x &= a_{-N} t^{-N} + a_{-N+1} t^{-N + 1} + \cdots \quad \text{and} \\
x-s &= b_{-N} t^{-N} + b_{-N+1} t^{-N + 1} + \cdots ,
\end{align*}
where $N$ is the smallest positive integer for which
$b_{-N}\neq 0$ (or $N=1$ if no such integer exists).
Then
$$\tau_{[s]}\psi(x) =
\mathbf{1}_{\FF_p[[t]]} (x-s) +
\left[\frac{(e^{2\pi i a_0/p} - 1)}{p^N (p-1)}\right]
\left[ e^{2\pi i c/p} + (p-1)e^{2\pi i d/p} \right],$$
where
\begin{align*}
c & = b_{-N+1} + b_{-N+2} + \cdots + b_{-1} \in \FF_p,
\quad \text{and}
\\
d & = b_{-N} + b_{-N+1} + \cdots + b_{-1} \in \FF_p.
\end{align*}

\end{example}

\begin{example}
\label{ex:Qpwave3}
Let $G=\QQ_3$, let $H=\ZZ_3$,
and let $A$ be multiplication-by-$1/3$, so that $|A|=3$.
Identify $\widehat{G}$ as $\QQ_3$ and $H^{\perp}$ as $\ZZ_3$,
as in Example~\ref{ex:Qp}, for $p=3$.
Let $\calD$ be the set of coset representatives for
$\widehat{G}/H^{\perp}$ defined in Example~\ref{ex:QpD}.
Define $V_{1/3}=(3\ZZ_3)\cup(2+3\ZZ_3)$, $V_{2/3}=1+3\ZZ_3$,
and $V_{\sigma}=\emptyset$
for all $\sigma\in \calD \setminus \{1/3,2/3\}$, and define
$T$ as in Theorem~\ref{thm:bb}.  That is, for all $\gamma\in 1+3\ZZ_3$,
$T(\gamma) = \gamma + 2/3$; and for all
$\gamma\in \ZZ_3\setminus (1+3\ZZ_3)$,
$T(\gamma) = \gamma + 1/3$.  Then for $n\geq 1$,
$$
\Lambda_n = \begin{cases}
1 +2\cdot 3 + 1\cdot 3^2 + 2\cdot 3^3 + \cdots
+ 1\cdot 3^{n-2} + 3^n\ZZ_3, &\text{if $n$ is even}; \\
2 +1\cdot 3 + 2\cdot 3^2 + 1\cdot 3^3 + \cdots
+ 1\cdot 3^{n-2} + 3^n\ZZ_3, &\text{if $n$ is odd}
\end{cases}
$$
Thus, $\Lambda_n$ is a closed ball of radius $1/3^{n}$
(and hence measure $1/3^n$) about a point which is in
$1+3\ZZ_3$ if $n$ is even, or in $2+3\ZZ_3$ if $n\geq 3$ is
odd, or at $0$ if $n=1$.
For $\Lambda$ and $\Omega$ as in equation~\eqref{eq:Omegadef2},
we see as before that
$\Lambda=\bigcup \Lambda_n \subset\ZZ_3$
is a disjoint union of countably many
balls, and that $\nu(\Lambda) = 1/(3-1)=1/2$.

Define $\psi=\check{\mathbf{1}}_{\Omega}$;
then by Theorem~\ref{thm:bb}, $\psi$ is a (single) wavelet.
And again, $\tau_{[s]}\psi$
must be constant on every set of the form $c + 3 \ZZ_3 \subset G$,
for any $[s]\in G/H$; but $\psi$ does not have compact support.
This time, when we compute $\tau_{[s]}\psi(x)$ by
equation~\eqref{eq:psicomp}, we get
\begin{align*}
\tau_{[s]}\psi(x) & =
\mathbf{1}_{\ZZ_3} (x-s)
\\
& +
\frac{\left[(x,2/3) - 1\right]}{3^{N}} 
(x-s,1+2\cdot 3 + 1\cdot 3^2 + \cdots+1\cdot 3^{(N-2)})
\\
& +
\frac{\left[(x,2/3) - 1\right]}{8\cdot 3^{N}} 
(x-s,1+2\cdot 3 + 1\cdot 3^2 + \cdots+2\cdot 3^{(N-1)})
\\
& +
\frac{\left[(x,1/3) - 1\right]}{8\cdot 3^{N-1}} 
(x-s,2+1\cdot 3 + 2\cdot 3^2 + \cdots+1\cdot 3^{(N-1)}),
\end{align*}
if $N$ is even, and
\begin{align*}
\tau_{[s]}\psi(x) & =
\mathbf{1}_{\ZZ_3} (x-s)
\\
& +
\frac{\left[(x,1/3) - 1\right]}{3^{N}} 
(x-s,2+1\cdot 3 + 2\cdot 3^2 + \cdots+1\cdot 3^{(N-2)})
\\
& +
\frac{\left[(x,1/3) - 1\right]}{8\cdot 3^{N}} 
(x-s,2+1\cdot 3 + 2\cdot 3^2 + \cdots+2\cdot 3^{(N-1)}),
\\
& +
\frac{\left[(x,2/3) - 1\right]}{8\cdot 3^{N-1}} 
(x-s,1+2\cdot 3 + 1\cdot 3^2 + \cdots+1\cdot 3^{(N-1)})
\end{align*}
if $N$ is odd, where $N$ is the smallest positive
integer such that $x-s\in 3^{-N}\ZZ_3$.
\end{example}

The author would like to thank John Benedetto and the referee
for their helpful comments on the exposition of this paper.

\bibliographystyle{amsalpha}

\bibliographystyle{plain}

\end{document}